\documentclass[12pt,a4j]{article}

\usepackage{amssymb}

\newenvironment{namelist}[1]{%
\begin{list}{}
{
\settowidth{\labelwidth}{#1}
\setlength{\leftmargin}{1.11\labelwidth}}
}{%
\end{list}}
\begin{document}
\begin{center}
{\bf\large Testing the equality of error
distributions from $k$ independent GARCH models
}
\end{center}
\begin{center}
{\sc {\sc S. Ajay
Chandra}}
\end{center}
\begin{center}
{\it Department of Mathematics and Statistics,\\
La Trobe University, Victoria, \\Australia}
\end{center}

\newcommand{\qed}{\hbox{\rule{6pt}{6pt}}}
\renewcommand{\abstractname}{{}}
\begin{abstract}
\noindent{\bf Abstract.}
 In this paper we study the problem of testing
the null hypothesis that errors from $k$
independent parametrically specified generalized
autoregressive conditional
 heteroskedasticity (GARCH) models have the
same distribution versus a general alternative.
First we establish the asymptotic validity of a
class of linear test statistics derived from the
$k$ residual-based empirical distribution
functions. A distinctive feature is that the
asymptotic distribution of the test statistics
involves terms depending on the distributions of
errors and the parameters of the models, and
weight functions providing the flexibility to
choose scores for investigating power
performance. A Monte Carlo study assesses the
asymptotic performance in terms of empirical size
and power of the three-sample test based on the
Wilcoxon and Van der Waerden score generating
functions in finite samples. The results
demonstrate that the two proposed tests have
overall reasonable size and their power is
particularly high when the assumption of Gaussian
errors is violated. As an illustrative example,
the tests are applied to daily individual stock
returns of the New York Stock Exchange data.
\end{abstract}
\section{Introduction}

Analysis of volatility in financial time series is certainly the
subject of considerable attention with huge literature having been
published. In the seminal papers by Engle (1982) and Bollerslev
(1986), generalized autoregressive conditional
 heteroskedasticity (GARCH) models have been
 proposed to capture special features of
financial volatilities.
 Since then, numerous variations and extensions of
GARCH models have been proposed to possibly
explain and model risk and uncertainty in pricing
derivative securities, in stochastic modelling of
 the term structure of interest rates, in applications
 related to fixed-income portfolio management, in
 asset pricing studies, and in the riskiness of financial returns which
provides a volatility measure that can be used in
financial decisions concerning risk analysis.
Several excellent surveys of the GARCH
methodology in finance are available, such as
Bollerslev et al. (1992), Engle (1995),
Gouri\'eroux (1997), Mikosch (2003) and Bauwens
{\it et al.} (2006).
\\
\indent For time series data, residuals must be
taken into account as they typically depend on
parameter estimates, and inference based on these
residuals, especially various diagnostic checks,
is a basic tool in the statistical analysis of
linear time series models (see Brockwell and
Davis (1994)). By contrast, asymptotic theory for
the residuals of nonlinear time series models has
been surveyed by Berkes and Horv\'ath (2002). For
a GARCH($p,q$) model, Berkes
 and Horv\'ath (2003) derived the asymptotic distribution
 of the empirical process of residuals and showed that,
unlike the residuals of autoregressive moving
average (ARMA) models, these residuals do not
behave in this context like asymptotically
independent random variables, and the asymptotic
distribution involves a term depending on the
parameters of
the model.\\
\indent The classical two-sample problem is one
of the central themes of nonparametric testing
theory. One of the problems most frequently
encountered in statistics is to test the
hypothesis of no difference between two
independent populations primarily on the basis of
samples drawn at random from these two
populations. Some of the earliest and most
classical tests of nonparametric nature for this
problem are Wilcoxon's test, the Mann and Whitney
test, the Mood and Brown test, Lehmann's test,
the Cram\'er-von Mises test and Van der Waerden's
test. Moreover, the classical limit theorem of
normalized two-sample linear test statistics
which generated much interest in this context is
the celebrated Chernoff$-$Savage (1958) theorem.
It is well known that the theorem is widely used
to study the asymptotic power and power
efficiency of the above two-sample tests. Further
refinements on their conditions of this theorem,
extensions and related results, are due to Durbin
(1973), Puri and Sen (1993) and
references therein. \\
\indent The natural extension of the two-sample
problem is the $k$-sample problem, where
observations are taken under a variety of
different and independent conditions. The
nonparametric test procedures which have been
developed for this $k$-sample problem require no
assumptions beyond continuous populations and
therefore are applicable under any circumstances.
The classical tests in this context are the
Kruskal-Wallis $H$ test, Terpestra's $k$-sample
test, the Mood and Brown $k$-sample test,
Kiefer's $k$-sample analogues of the
Kolomogorov-Smirnov test and the Cram\'er-von
Mises $k$-sample test. To this end, it is of
interest to state that Puri (1964) generalized
the situation covered by the Chernoff$-$Savage
(1958) theorem to the
$k$-sample problem.\\
\indent If GARCH errors were observable, the
problem that we consider here would be the
classical $k$-sample problem studied by Puri
(1964). In our context, we do not observe these
errors, but assume that well-behaved estimators
of the parameters of the model are available.
Hence, our test procedure can be thought of as an
extension of the $k$-sample problem. More
specifically, we are concerned with testing the
null hypothesis that errors from $k$ independent
parametrically specified GARCH models have the
same distribution versus a general alternative in
the spirit of Chernoff and Savage (1958), Puri
(1964), and Berkes and Horv\'ath (2003). In
contrast with the independent, identically
distributed or ARMA setting, this study
highlights some interesting features
of $k$ GARCH residual-based test statistics. \\
\indent Potential applications of the $k$-sample
test are to be found especially in studies of the
behavior of speculative prices, such as stock
prices or exchange rates, usually in view of
testing market efficiency. One important problem,
for example, the stock return of a company is
defined as the error from a GARCH model, and the
researcher is often interested in comparing the
distributions of stock return of companies from
$k$ independent groups. Another related problem
in this context is that the researcher may be
interested in comparing the distribution
functions of standardized real variables like
exports or output growth rates with data from $k$
independent companies. In other areas of
financial markets it is often of interest to test
whether $k$ observable variables belong to the
same location-scale family, which is also a
special case of the test that we study. In all
these situations, the usual approach to test for
the equality of the distribution functions is to
test the equality of just some moments to propose
parametric models for the errors and then test
whether the parameters estimated are equal.
Instead, we propose to compare the entire
distribution functions without assuming any
parametric form for them.  \\
\indent The objective of this paper is to study
the asymptotic behavior of $k$ GARCH
residual-based linear test statistics. The rest
of the paper is organized as follows. Section 2
introduces the construction of $k$ GARCH
residual-based empirical distribution functions
and proposes linear test statistics pertaining to
these residual-based empirical distribution
functions. In Section 3, we establish the
asymptotic validity of the test. Section 4
reports the results in terms of empirical size
and power of a Monte Carlo study for
 validating the three-sample test based on the Wilcoxon
and Van der Waerden score generating functions
for finite sample sizes. As an example, the two
tests are applied to daily individual stock
returns of the New York Stock Exchange data. The
proof of the result in Section 3 is provided in
Section 5.
\section{ $k$ GARCH residual-based linear test statistics} In this section, we propose a family of linear
test statistics pertaining to empirical processes
of residuals in order to test the null hypothesis
that errors from $k$ parametrically specified
GARCH models have the same distribution against a
general alternative. We shall formulate the
$k$-sample problem as follows. Let us consider
the $k$ independent random samples generated from
the GARCH($p_j,q_j$) models given by
\begin{eqnarray}
 \left.
\begin{array}{lll}
X_{j,t}&=&\sigma_{j,t}\varepsilon_{j,t},\\
 \sigma_{j,t}^2&=&
\omega_{0j}+\sum\limits_{i=1}^{p_j} \alpha_{0j}^i
X_{j,t-i}^2+\sum\limits_{i'=1}^{q_j}
\beta_{0j}^{i'} \sigma_{j,t-i'}^2,\quad 1\le t\le
n_j,\quad 1\le j\le k,
\end{array}
\right.
\end{eqnarray}
where the $\varepsilon_{j,t}$ are independent and
identically distributed random variables such
that $E(\varepsilon_{j,t}^2)=1$, $\omega_{0j}>0$,
$\alpha_{0j}^i\ge 0$, $1\le i\le p_j$,
$\beta_{0j}^{i'}\ge 0$, $1\le i'\le q_j$, and the
$\varepsilon_{j,t}$ is independent of
$X_{j,s},s<t$. Henceforth, it is tacitly assumed
that $\alpha_{0j}^{p_j}>0$ when $p_j\ge 1$, and
$\beta_{0j}^{q_j}>0$ when $q_j\ge 1$.
 \\
 \indent In this paper, we are primarily concerned
 with the $k$-sample problem of testing
\begin{eqnarray}
 \left.
\begin{array}{l}
H_0:F_1(x)=\cdots=F_k(x)\;\mbox{for all
$x$}\\
\hspace*{-3.4cm}\mbox{against}  \\
H_A:\mbox{$F_i(x)\ne F_j(x)$ for at least some
$x$, and $i\ne j$},
\end{array}\right.
\end{eqnarray}
where $F_j(\cdot)$ is the distribution function
of $\{\varepsilon_{j,t}\}$, which is assumed to
be absolutely continuous with respect to the
Lebesgue measure, but unspecified. Henceforth, we
assume that $f_j(x)=F_j'(x)$ exists and is
defined over
$(-\infty,\infty)$.  \\
 \indent We first proceed to describe the quasi-maximum
likelihood (QML) estimation of model (1). The
vector of parameters is
$\theta_j=(\theta_{1,j},\ldots,
\theta_{j,p_j+q_j+1})^T =(\omega_j,
\alpha_j^1,\ldots,\alpha_j^{p_j},\beta_j^1,
\ldots,\beta_{j}^{q_j})^T$ which belongs to a
compact parameter space $\Theta_j\subset
(0,\infty)^j\times [0,\infty)^{p_j+q_j}$. The
true vector of parameters is unknown and is
denoted by $\theta_{0j}=(\omega_{0j},
\alpha_{0j}^1,\ldots,\alpha_{0j}^{p_j},\beta_{0j}^1,
\ldots,\beta_{0j}^{q_j})^T$.\\
\indent Suppose that an observed stretch
$X_{j,1},\ldots,X_{j,{n_j}}$ from $\{X_{j,t}\}$
is available. Note that if
$\{\varepsilon_{j,t}\}$ is Gaussian, the
quasi-likelihood function with respect to initial
values $X_{j,0},\ldots,X_{j,1-p_j}$,
$\tilde{\sigma}_{j,0}^2,\ldots,
\tilde{\sigma}_{j,1-q_j}^2$, is given by
\[\mathbb{L}_{n_j}(\theta_j)=\sum_{t=1}^{n_j} \frac{1}
{\sqrt{2 \pi \tilde{\sigma}_{j,t}^2
}}\exp\biggl(-\frac{X_{j,t}^2}{2
\tilde{\sigma}_{j,t}^2}\biggl),
\]
where the $\tilde{\sigma}_{j,t}^2$, $t\ge 1$ are
defined recursively by
\[
\tilde{\sigma}_{j,t}^2=\tilde{\sigma}_t^2(\theta_j)=\omega_j+\sum
\limits_{i=1}^{p_j} \alpha_j^i
X_{j,t-i}^2+\sum\limits_{i'=1}^{q_j} \beta_j^{i'}
\tilde{\sigma}_{j,t-i'}^2,\qquad 1\le j\le k.
\]
As an example, one can choose the initial values
as
$X_{j,0}^2=\cdots=X_{j,1-p_j}^2=\tilde{\sigma}_{j,0}^2=\cdots=
\tilde{\sigma}_{j,1-q_j}^2\equiv\omega_j$ or
$X_{j,0}^2=\cdots=X_{j,1-p_j}^2=\tilde{\sigma}_{j,0}^2=\cdots=
\tilde{\sigma}_{j,1-q_j}^2=X_{j,1}^2$. \\
\indent We can now define the QML estimators of
$\theta_j$ by
\[
\hat {\theta}_{j,n_j}= \arg \max_{\theta_j\in \Theta_j}
  \mathbb{L}_{n_j}(\theta_j)=
\arg \min_{\theta_j\in \Theta_j}
  \tilde{\mathcal{I}}_{n_j}(\theta_j),
\]
where
\[
\tilde{\mathcal{I}}_{n_j}(\theta_j)=\frac{1}{n_j}
\sum_{t=1}^{n_j}
 \tilde{l}_t(\theta_j),\qquad \tilde{l}_t(\theta_j)
 =\log \tilde{\sigma}_{j,t}^2
 +\frac{X_{j,t}^2}{\tilde{\sigma}_{j,t}^2},\qquad
 1\le j\le k.
\]
For $\hat {\theta}_{j,n_j}$, it is assumed that
\begin{equation}
\|\hat
{\theta}_{j,n_j}-\theta_{0j}\|=\mathcal{O}_p(n_j^{-1/2}),
\qquad 1\le j\le k,
\end{equation}
where $\|\cdot\|$ denotes the Euclidean norm. The
validity of (3) is established by Francq and
Zako\"ian (2004) based on the conditions of
Assumption 2 given below. Conditions (3) are also
typically satisfied by the QML estimators of
Straumann and Mikosch (2006). Henceforth, the
empirical residuals are given by
$$\hat{\varepsilon}_{j,t}=X_{j,t}/
\tilde{\sigma}_t(\hat{\theta}_{j,n_j}),\qquad
1\le j\le k.$$ \indent For (2), we first collect
some basic tools and then describe our approach
in the spirit of Chernoff and Savage (1958), and
Puri (1964). Write $N=\sum_{j=1}^k n_j$ and
$\lambda_{jN}=n_j/N$, $1\le j\le k$. In the
following, we assume that the inequalities
$0<\lambda_0\le
\lambda_{1N},\ldots,\lambda_{kN}\le
1-\lambda_0<1$ for some $\lambda_0\le 1/k$.
Define by
\[
H_N(x)=\sum_{j=1}^k
\lambda_{jN}F_j(x)
\]
 the combined cumulative distribution function.
 Write
$F_{j,n_j}(x)=n_j^{-1}\sum_{t=1}^{n_j}
[I(\varepsilon_{j,t}\le x)]$ and
$\hat{F}_{j,n_j}(x)=n_j^{-1}\sum_{t=1}^{n_j}
[I(\hat{\varepsilon}_{j,t}\le x)]$, where
$I(\Omega)$ is the indicator function of the
event $\Omega$. Then the empirical distribution
function is
\[\mathcal{H}_N(x)
=\sum_{j=1}^k\lambda_{jN}F_{j,n_j}(x)
\]
and analogously,
\begin{equation}
\hat{\mathcal{H}}_N(x)=\sum_{j=1}^k\lambda_{jN}
\hat{F}_{j,n_j}(x).
\end{equation}
Set $\hat{B}_{j,n_j}(x)
=n_j^{1/2}(\hat{F}_{j,n_j}(x)-F_j(x))$. Then by
virtue of Berkes and Horv\'ath (2003), it follows
that
\begin{equation}
\hat{B}_{j,n_j}(x)=\mathcal{E}_{j,n_j}(x)
+\mathcal{A}_jxf_j(x)+ \xi_{j,n_j}(x),
\end{equation}
where $\sup_x|\xi_{j,n_j}(x)|=o_p(1)$,
\[
\mathcal{E}_{j,n_j}(x)=n_j^{-1/2}\sum_{t=1}^{n_j}
[I(\varepsilon_{j,t}\le x)-F_j(x)],\quad
\mathcal{A}_j=\sum_{l=1}^{p_j+q_j+1}n_j^{1/2}
(\hat{\theta}_{j,n_j}^l-\theta_{0j}^l)\tau_{j,l}
\]
with
\[\tau_{j,1}=E[1/2\tilde{\sigma}_t^2(\theta_{0j})],
\quad \tau_{j,l}=E[
X_{j,t-l}^2/2\tilde{\sigma}_t^2(\theta_{0j})],
\quad 2\le l\le p_j+1, \] and
$\tau_{j,p_j+1+l'}=E[\tilde{\sigma}_{t-l'}^2(\theta_{0j})
/2\tilde{\sigma}_t^2(\theta_{0j})]$, $1\le l'\le
q_j$, $1\le j\le k$. Hence, by analogy with (5),
the asymptotic representation of (4) becomes
\begin{equation}
\hat{\mathcal{H}}_N(x)=\mathcal{H}_N(x)+\sum_{j=1}^k
n_j^{-1/2}\lambda_{jN}\mathcal{A}_jx f_j(x)
+o_p(N^{-1/2}).
\end{equation}
Decomposition
(6) is basic and plays an important role in the sequel.\\
\indent Define $\hat{S}_{iN}^{(j)}=1$, if the
$i$th smallest of $N=\sum_{j=1}^k n_j$ empirical
residuals is from $\{\hat{\varepsilon}_{j,t}\}$,
and otherwise define $\hat{S}_{iN}^{(j)}=0$,
$1\le i\le N$, $1\le j\le k$. Then, for (2), we
shall consider a family of linear test statistics
of the form
\[
\hat{T}_{jN}=\frac{1}{n_j}\sum_{i=1}^N
E_{iN}\hat{S}_{iN}^{(j)},\quad 1\le j\le k,
\]
 where the $E_{iN}$ are given constants
called weights or scores. The definition of
$\hat{T}_{jN}$ is the one traditionally used. We
shall, however, use the representation given by
\begin{equation}
\hat{T}_{jN}=\int
J\biggl(\frac{N}{N+1}\hat{\mathcal{H}}_N(x)\biggl)d
\hat{F}_{j,n_j}(x),\quad 1\le j\le k,
\end{equation}
where $J(u)$, $0<u<1$, is a continuous
score-generating function. Note that
$E_{iN}=J(i/(N+1))$, $1\le i\le N$ are functions
of the ranks $i$ ($=1,\ldots,N$) and are
explicity known. Some typical examples of $J$
given in Puri and Sen (1993) are as follows:
\begin{itemize}
 \item[(i)] Wilcoxon's $k$-sample test with $J(u)=u$, $0<u<1$,
\item[(ii)] Van der Waerden's $k$-sample test
with $J(u)=\Phi^{-1}(u)$, $0<u<1$, where
$\Phi(x)=(2\pi)^{-1/2}\int_{-\infty}^x
e^{-t^2/2}dt$, \item[(iii)] Mood's $k$-sample
test with $ J(u)=(u-\frac12)^2$, $0<u<1$,
\item[(iv)] Klotz's normal $k$-sample test with
$J(u)=(\Phi^{-1}(u))^2$, $0<u<1$.
\end{itemize}

\indent In the following, $K$ will denote a
generic constant taking many different values
$K>0$ which may depend on $J$ but will not depend
on $F_j(\cdot)$, $n_j$ and $N$ for all $1\le j\le
k$.

\section{Asymptotic properties of $\hat{T}_{jN}$}
In this section, our primary object is to show
that (7) has an asymptotically normal
distribution. For this purpose, let
$\{\Delta_{0j,t},1\le j\le k\}$ be the $(p_j+ q_j
)\times(p_j+ q_j)$ matrices defined by
\[
\Delta_{0j,t}=\left(
\begin{array}{cccccccc}
\alpha_{0j}^1\varepsilon_{j,t}^2&\cdots&\alpha_{0j}^{p_j-1}
\varepsilon_{j,t}^2&\alpha_{0j}^{p_j}\varepsilon_{j,t}^2&
\beta_{0j}^1\varepsilon_{j,t}^2&\cdots&\beta_{0j}^{q_j-1}
\varepsilon_{j,t}^2&\beta_{0j}^{q_j}\varepsilon_{j,t}^2\\
& I_{p_j-1}&&&&0_{(p_j-1  )\times(q_j+ 1) }&\\
\alpha_{0j}^1&\cdots&\alpha_j^{p_j-1}&\alpha_{0j}^{p_j}
&\beta_{0j}^1&\cdots&\beta_{0j}^{q_j-1}&\beta_{0j}^{q_j}
\\
&&0_{(q_j-1)\times p_j}&&& I_{q_j-1}&&0_{(q_j-
1)\times 1}
\end{array}
\right).
\]
Assuming that
\begin{equation}
E(\log^+\|\Delta_{0j,1}\|)\le
E\|\Delta_{0j,1}\|<\infty,
\end{equation}
the top Lyapunov exponent is defined by
$\gamma(\Delta_{0j})\equiv \inf _{t\ge 1}t^{-1}
E(\log \|\Delta_{0j,1}\Delta_{0j,2}\cdots
\Delta_{0j,t} \|)$, where
$\Delta_{0j}=\{\Delta_{0j,t},1\le j\le k\}$. In
particular, one can readily check that if
$\{\varepsilon_{j,t}\}$ is Gaussian, (8) holds.
Bougerol and Picard (1992a,b) showed that if (8)
holds, a general GARCH($p_j,q_j$) process has a
unique non-anticipative strictly stationary
solution if and only if $\gamma(\Delta_{0j})<0$,
$1\le j\le k$. \\
\indent To establish the asymptotic properties of
(7), we impose the following
regularity conditions. \\\\
\noindent{\bf Assumption 1}
\begin{itemize}
 \item [(A.1)] $J(u)$ is not constant and has a continuous derivative $J'(u)$ on (0,1).
\item[(A.2)]$|J(u)|\le
K[u(1-u)]^{-\frac12+\delta}$ and $|J'(u)|\le
K[u(1-u)]^{-\frac32+\delta}$ for some $\delta>0$.
\item[(A.3)] $xf_j(x)$ and $xf_j'(x)$ are
uniformly bounded continuous, and integrable
functions on $(-\infty,\infty)$.  \item[(A.4)]
There exist constants $c_j>0$ such that
$F_j(x)\ge c_j\{xf_j(x)\}$ for all $x>0$.
 \end{itemize}

\indent A few remarks concerning the necessity of
these conditions are in order. Assumptions (A.1)
and (A.2) are basic conditions in our context. As
noted by Chernoff and Savage (1958), typically
(A.2) has two important functions: (i) it limits
the growth of the function $J$ and (ii) it
supplies certain smoothness properties. Both
conditions can be easily verifiable in the
preceding examples given by $J$. Assumption (A.3)
is basic and necessary for studying residual
empirical processes and establishing the
convergence result of (7). This condition was
also made for empirical processes pertaining to
linear regression residuals by Bai (1996).
Assumption (A.4) is virtually imposed in dealing
with the convergence of higher order terms of
(7). Finally, it is worth noting that conditions
(A.1)$-$(A.4) are typically satisfied by several
error distributions such as, normal, Student's
$t$, logistic, double exponential, gamma and
Laplace.
\\
\indent To validate (3), we require the following additional
regularity conditions, which can be
found in Francq and Zako\"ian (2004).  \\

 \noindent{\bf Assumption 2}
\begin{itemize}
\item[(B.1)] $\theta_{0j}\in \tilde{\Theta}_j$,
where $\tilde{\Theta}_j$ denotes the interior of
the compact parameter space $\Theta_j$.
\item[(B.2)] $\gamma(\Delta_{0j})<0$ and
$\sum_{i'=1}^{q_j} \beta_{0j}^{i'}<1$ for all
$\theta_j\in \Theta_j$.
 \item[(B.3)] $\varepsilon_{j,t}^2$ has a
 non-degenerate distribution with
$E(\varepsilon_{j,t}^2)=1$.
 \item[(B.4)]
$\kappa_j\equiv E(\varepsilon_{j,t}^4)<\infty$.
\item[(B.5)] If $q_j>0$, $A_{\theta_{0j}}(z)$ and
$B_{\theta_{0j}}(z)$ have no common root,
$A_{\theta_{0j}}(1)\not= 0$, and
$\alpha_{0j}^{p_j} +\beta_{0j}^{q_j}\not=0$,
where $A_{\theta_{0j}}(z)=\sum_{i=1}^{p_j}
\alpha_j^iz^i$ and $B_{\theta_{0j}}(z)
=1-\sum_{i'=1}^{q_j} \beta_j^{i'}z^{i'}$.
Conventionally, $A_{\theta_{0j}}(z)=0$
 if $p_j=0$ and $B_{\theta_{0j}}(z)=1$ if
 $q_j=0$.
\end{itemize}

\indent We now justify that conditions
(B.1)$-$(B.5) are necessary for the
 model under consideration. These conditions were
 essentially made by Francq and Zako\"ian (2004)
 for the validity of (3). We first note that the
 compactness of $\Theta_j$ is always assumed.\\
 \indent Assumption (B.1) is typically necessary to
 obtain the asymptotic normality of the QML estimators
 $\hat{\theta}_{j,n_j}$, $1\le j\le k$. In the case of
$\alpha_{0j}\equiv\alpha_{0j}^1=0$, the limit
distribution of
$\sqrt{n_j}(\hat{\alpha}_j-\alpha_{0j})$ is
non-normal over $[0,\infty)$. Assumption (B.2) is
a sufficient condition for the stationarity and
ergodicity of model (1). This condition implies
that the roots of $B_{\theta_{j}}(z)$ are outside
the unit disc. Moreover, if
$\gamma(\Delta_{0j})<0$, there exists $s>0$ such
that $E(\sigma_{j,t}^{2s})<\infty$ and
$E(X_{j,t}^{2s})<\infty$. Assumption (B.3) is
made for model identification is not restrictive
provided $E(\varepsilon_{j,t}^{2})<\infty$. This
moment condition is clearly necessary to
establish the asymptotic normality of the
Gaussian QML estimator as in Berkes and
Horv\'{a}th (2003). The existence of a
fourth-order moment given by (B.4) is a
strengthening of (B.3) required for the
finiteness of the variance of the score vector
$\partial \tilde{l}_t(\theta_{0j})/\partial
\theta_j$. Note also that this condition does not
imply the existence of a second-order moment for
the observed process $\{X_{j,t}\}$. It is often
the case that the existence of the second-order
moments is found to
be inappropriate for financial applications.\\
\indent Finally, the assumption that the polynomials whose common
roots uniquely identify $\theta_j$ was also made by Berkes {\it et.
al} (2003). This condition is typically satisfied when $p_j>1$ and
$q_j>1$. If $p_j=1$ and $\alpha_{0j}\not =0$, the unique root of
$A_{\theta_{0j}}(z)=0$ and $B_{\theta_{0j}}(z)\not=0$. If $q_j=1$
and $\beta_{0j}\equiv \beta_{0j}^1\not =0$, the unique root of
$B_{\theta_{0j}}(z)=1/\beta_{0j}>0$, and because $\alpha_{0j}>0$
produces $A_{\theta_{0j}}(1/\beta_{0j})\not=0$. Moreover, it can be
noted that (B.5) implies that $\theta_{0j}$ does not necessarily
have to belong to the interior of $\Theta_j$. This is essentially
important when dealing with situations of over-specification. When a
GARCH($p_j,q_j$) is fitted, one can show that an ARCH($p_j$) model
can be estimated consistently. In a general sense, either $p_j$ or
$q_j$ can be over-specified, but not both of them. Indeed, it is
required that $\alpha_{0j}^i>0$ for some $i$ when $p_j>0$. If this
assumption is dropped, the model solution would simply reduce to an
i.i.d. white noise of the form
$\sigma_{j,t}^2=\sigma_j^2(1-\beta_{0j})
+\beta_{0j}\sigma_{j,t-1}^2$, where
$\sigma_j^2=\omega_{0j}/(1-\beta_{0j})$.\\
\indent In order to state the main result, we
shall introduce the following notation:
\[
\mathcal{U}(\theta_{0j})= E
\biggl[\frac{1}{\sigma_t^4(\theta_{0j})}
\frac{\partial \sigma_t^2(\theta_{0j})}{\partial
\theta_j}\frac{\partial \sigma_t^2
(\theta_{0j})}{\partial \theta_j^T}\biggl],\quad
u_{t}(\theta_j)=
\frac{1}{\sigma_t^2(\theta_j)}\frac{\partial
\sigma_t^2(\theta_{0j})}{\partial \theta_j},\quad
1\le j\le k. \]
 By virtue of (B.4), it is seen that the $i$th element of each $\hat
{\theta}_{j,n_j}$, $1\le j\le k$ admits the
asymptotic representation,
\[
\hat
{\theta}_{j,n_j}^i-\theta_{0j}^i=\frac{1}{n_j}
\sum_{t=1}^{n_j}
Z_{t}^i(\theta_j)(\varepsilon_{j,t}^2-1)
+o_p(n_j^{-1/2}),
\]
where $Z_{t}^i(\theta_j)$ is the $i$th element of
$[\mathcal{U}(\theta_{0j})]^{-1}u_{t}(\theta_j)$,
$1\le i\le p_j+q_j+1$. As shown by Francq and
Zako\"ian (2004), $\mathcal{U}(\theta_{0j})$ is
positive definite for all $1\le j\le k$. These
considerations motivate the following result,
whose proof is
relegated to Section 5.\\   \\
 \noindent{\bf Theorem
1.} {\it Suppose that Assumptions 1 and 2 hold
and that, in addition, $\{\hat{\theta}_{j,n_j}\}$
is a sequence of QML estimators typically
satisfying (3). Then, as $N\to\infty$,
\[
N^{1/2}\Sigma_N^{-1/2} s_N\stackrel{d}
{\longrightarrow}\mathcal{N}(0,I_k),
\]
where $I_k$ is the $k\times k$ identity matrix,
$\Sigma_N$ is the $k\times k$ positive definite
dispersion matrix whose entries are given by (15)
and (16), and $s_N= (\hat{T}_{jN}-\mu_{jN})
_{1\le j\le k}$ with
$\mu_{jN}=\int J(H_N)dF_j(x)$.} \\\\
\noindent {\bf Remark 1.} If $J(\cdot)$ and
$\Sigma_N$ were known, an immediate consequence
of Theorem 1 is that the quadratic statistic
$\mathcal{L}_N=N s_N^T\Sigma_N^{-1}s_N$ has an
approximate $\chi^2(k)$ distribution with $k$
degrees of freedom under $H_0$ (cf. Theorem 2.8
in Seber (1977)). Unfortunately, the covariance
structure of $\Sigma_N$, in general, depends on
the unspecified distribution function
$F_j(\cdot)$, the unknown parameter vector
$\theta_{0j}$ and some expectations. Thus, it is
not possible to perform a consistent test based
on $\mathcal{L}_N$. Replacing $\Sigma_N$ by a
consistent estimator $\hat{\Sigma}_N$ (for
details, see Section 4), we can effectively
estimate $\mathcal{L}_N$ by
$\hat{\mathcal{L}}_N=Ns_N^T\hat{\Sigma}_N^{-1}s_N$.
Writing $\tilde{s}_N=N^{1/2}
(\hat{T}_{jN}-\mu_{jN}) _{1\le j\le k}$, we have
$\hat{\mathcal{L}}_N/\mathcal{L}_N= \tilde{s}_N^T
\hat{\Sigma}_N^{-1}\tilde{s}_N/\tilde{s}_N^T
\Sigma_N^{-1}\tilde{s}_N$, and using Lemma 1
given in Section 5, it
 follows that $\mbox{ch}_k(\hat{\Sigma}_N\Sigma_N^{-1})\le
(\hat{\mathcal{L}}_N/\mathcal{L}_N)\le
\mbox{ch}_1(\hat{\Sigma}_N\Sigma_N^{-1})$, where
$\mbox{ch}_j(\Lambda)$ is the $j$th
characteristic root of $\Lambda$. Moreover, by
the ergodic theorem we have
$\hat{\Sigma}_N\Sigma_N^{-1}\stackrel{p}{\rightarrow}
 I_k$, which implies $\mbox{ch}_1(\hat{\Sigma}_N\Sigma_N^{-1})
  \stackrel{p}{\rightarrow} 1$ and $\mbox{ch}_k(\hat{\Sigma}_N\Sigma_N^{-1})
  \stackrel{p}{\rightarrow} 1$.
Observing that $\hat{\mathcal{L}}_N/\mathcal{L}_N
 \stackrel{p}{\rightarrow} 1$, and writing $\hat{\mathcal{L}}_N
 =\mathcal{L}_N \times
 (\hat{\mathcal{L}}_N/\mathcal{L}_N)$ we may
 conclude from Slutsky's theorem that $\hat{\mathcal{L}}_N
\stackrel{d}{\rightarrow}\chi^2(k)$ under $H_0$,
as was to be proved.
 \section{Simulation and empirical studies}
In this section we study the finite sample
performance of the proposed test procedure by
means of a simple numerical experiment and an
empirical example. The ideal way to carry out the
former case would be first to generate data from
some specific GARCH model, and then estimate a
GARCH model either correctly specified or not and
check the asymptotic behavior of
$\hat{\mathcal{L}}_N$ in terms of empirical size and power.  \\
\indent For simplicity and clarity, we shall
consider three-independent random samples
generated from the GARCH(1,1) model
\begin{eqnarray}
 \left.
\begin{array}{lll}
X_{j,t}&=&\sigma_{t}(\theta_j)\varepsilon_{j,t},\quad
 \sigma_{t}^2(\theta_j)=
\omega_{j}+ \alpha_{j} X_{j,t-1}^2+ \beta_{j}
\sigma_{t-1}^2(\theta_j),\quad 1\le t\le
n_j,\quad 1\le j\le 3,
\end{array}
\right.
\end{eqnarray}
where the $\varepsilon_{j,t}$ are independent and
identically distributed random variables such
that $E(\varepsilon_{j,t}^2)=1$,
$\theta_j=(\omega_j,\alpha_j,\beta_j)^T$,
$\omega_{j}>0$, $\alpha_{j}\ge 0$, $\beta_{j}\ge
0$ are unknown parameters, and the
$\varepsilon_{j,t}$ are independent of
$X_{j,s},s<t$. Note that model (9) is the most
commonly used in the literature, and enjoy
substantial application
in the finance setting.\\
\indent In the following, we are concerned with
the three-sample problem of testing
\[
H_0:F_1(\cdot)=F_2(\cdot) =F_3(\cdot)\qquad
\mbox{against}  \qquad H_A:F_1(\cdot)\ne
F_2(\cdot)\ne F_3(\cdot),
\]
where $F_j(\cdot)$ is an absolutely continuous
distribution function of $\{\varepsilon_{j,t}\}$,
but unspecified. For testing $H_0$, we propose to
use the statistic
$\hat{\mathcal{L}}_N=Ns_N^T\hat{\Sigma}_N^{-1}s_N$,
which has an approximate $\chi_r^2(3)$
distribution with $3$ degrees of freedom and
$0<r<1$ is the preassigned level of significance.
\\
 \indent We now describe our
 goodness-of-fit test using a smoothed
 bootstrap procedure. To this end, note that the
 asymptotic distribution of $\hat{T}_{jN}$
 depends crucially on the assumption of continuity
 and hence bootstrap samples must be
 generated from continuous distributions.
The following steps provide an explicit
description of the bootstrap test procedure based
on $\hat{\mathcal{L}}_N$:
 \begin{namelist}{4.}
    \item [1.] Having observed $X_{j,1},\ldots,X_{j,n_j}$, obtain an
    estimate $\hat{\theta}_{j,n_j}=(\hat{\omega}_j,
    \hat{\alpha}_j,\hat{\beta}_j)^T$ of $\theta_j$ using the QML
    method described in Section 2.
    \item [2.] Generate $B$ independent sequences
    of i.i.d. standard normal random variables
with replacement, each of length
    $n_j+n_0$, where $n_0$ is the length of warm-up sequence to reduce the effect of initial
conditions. Then define each of
    the $B$ sequences by $\varepsilon_{j,-n_0+1}^*,\ldots,
    \varepsilon_{j,0}^*,\varepsilon_{j,1}^*,\ldots,
    \varepsilon_{j,n_j}^*$.
    \item [3.] Generate $B$ bootstrap GARCH(1,1) independent
realizations $X_{j,1}^*,\ldots,X_{j,n_j}^*$ with
replacement, where the $X_{j,t}^*$, by analogy
with (9), satisfy
    \[
      X_{j,t}^*=\sigma_{t}^*(\hat{\theta}_{j,n_j})
      \varepsilon_{j,t}^*,\quad
      \sigma_{t}^{* 2}(\hat{\theta}_{j,n_j})=
      \hat{\omega}_j+\hat{\alpha}_j
      X_{j,t-1}^{* 2}+\hat{\beta}_j
      \sigma_{t-1}^{* 2}(\hat{\theta}_{j,n_j}).
    \]
    Note that $\{X_{j,t}^*\}$ is a smooth
    bootstrap version of the sample
    $\{X_{j,t}\}$.
    \item [4.] For each of the $B$ samples
    $X_{j,1}^*,\ldots,X_{j,n_j}^*$, obtain an
    estimate $\hat{\theta}_{j,n_j}^*
=(\hat{\omega}_j^*$, $\hat{\alpha}_j^*$,
    $\hat{\beta}_j^*)^T$ of $\hat{\theta}_{j,n_j}$ and construct bootstrap
    empirical residuals
    \[
    \hat{\varepsilon}_{j,t}^*=X_{j,t}^*\biggl/
\sqrt{\hat{\omega}_j^*
    +\hat{\alpha}_j^*X_{j,t-1}^{* 2}+\hat{\beta}_j^*
    \sigma_{t-1}^{* 2}(\hat{\theta}_{j,n_j}^*)},
\quad t=2,\ldots,n_j,\quad 1\le j\le 3.
    \]
    \item [5.] For the score generating functions $J(u)=u$
(Wilcoxon) and $J(u)=\Phi^{-1}(u)$ (Van der
Waerden), evaluate the following integral by a
rectangular numerical integration with $m$ terms:
\[
\hat{T}_{jN}^*=\int
J\biggl(\frac{N}{N+1}\hat{\mathcal{H}}_N^*(x)\biggl)d
\hat{F}_{j,n_j}^*(x),\quad 1\le j\le 3,
\]
where $\hat{F}_{j,n_j}^*(\cdot)$ denotes the
empirical distribution function constructed from
$\{\hat{\varepsilon}_{j,t}^*\}$ and
$\hat{\mathcal{H}}_N^{*} (\cdot)$ is the
bootstrap version of (6). Then, for each of the
$B$ residuals $\{\hat{\varepsilon}_{j,t}^*\}$,
calculate $\hat{\mathcal{L}}_N^*=N s_N^{*
T}\hat{\Sigma}_N^{* -1}s_N^*$, where
$s_N^*=(\hat{T}_{jN}^*-\mu_{jN}) _{1\le j\le 3}$
and $\hat{\Sigma}_N^*$ is a resampled version of
$\hat{\Sigma}_N$.
    \item [6.] Finally, repeat step 5 $B$ times and
then reject $H_0$ with significance level $r$ if
the $p$-value $\hat{r}=P(\hat{\mathcal{L}}_{N}^*>
\mathcal{L}_N^*)<r$, where $\mathcal{L}_N^*$ is
the $1-r$ sample quantile from
$\{\hat{\mathcal{L}}_{N,b}^*\}_{b=1}^B$. Here $B$
is chosen to be a sufficiently large integer.
 \end{namelist}

\indent In what follows we test the null
hypothesis that the zero-mean unit-variance
errors have the same distribution function at the
5\% significance level. For this purpose, we
shall consider two data generating processes
(DGPs):
   \begin{eqnarray*}
 \left.
\begin{array}{lllrcl}
  &&&X_{j,t}&=&\sigma_{t}
(\theta_j)\varepsilon_{j,t},\\
\vspace*{0cm}   \\
 \mbox{DGP 1:}&&&
 \sigma_{t}^2(\theta_j)&=&
0.1+ 0.1\, X_{j,t-1}^2+ 0.1\,
\sigma_{t-1}^2(\theta_j), \\
\vspace*{0cm}   \\
    \mbox{DGP 2:}&&&
 \sigma_{t}^2(\theta_j)&=&
0.5+ 0.4\, X_{j,t-1}^2+ 0.4 \,
\sigma_{t-1}^2(\theta_j),\quad 1\le j\le 3,
\end{array}
\right.
\end{eqnarray*}
where the $\varepsilon_{1,t}$ are i.i.d. random
variables with an $\mathcal{N}(0,1)$
distribution, the $\varepsilon_{2,t}$ are i.i.d.
random variables with mixture distribution
$(1-\varphi) \mathcal{N}(0,1)+\varphi
\mathcal{N}(2,1)$, $0\le \varphi\le 1$ and the
$\varepsilon_{3,t}$ are i.i.d. random variables
with Student's $t$ distribution having
$\varphi^{-1}$ degrees of freedom. The values of
$\varphi$ that we consider are $\varphi\in \{0,
1/9, 1/5, 1/3\}$. Note that if $\varphi=0$, the
errors $\varepsilon_{2,t}$ and
$\varepsilon_{3,t}$ are generated from a standard
normal distribution. The choice of $\varphi$
values, in principle, indicates that the last two
error processes have a leptokurtic distribution
whose tails are heavier than the ones of a normal
distribution. Observe that $H_0$ holds true if
and only if $\varphi=0$. We also notice that the
parameter $\varphi$ represents the departure from
$\mathcal{N}(0,1)$ in the sense that the larger
the value of $\varphi$, the larger the deviation
from the null model. Here, the distributions of
interest are re-scaled such that they have the
required zero mean and unit variance.
 \\
\indent We generate repeated trials of lengths
$n_1=n_2=n_3\in\{100,300,500\}$ from DGP1 and
DGP2, and compute the empirical size and power of
the 3$-$sample bootstrap Wilcoxon (W) and Van der
Waerden (VdW) tests at the 5\% nominal level
based on the steps 1$-$6 for each trial. The
number of Monte Carlo trials is 10000 with
$B=1000$ bootstrap replications each. Each
configuration of parameters was estimated by the
QML method.\\
\indent Table 1 reports the empirical proportion
of rejections of $H_0$ for the  W and VdW tests
based on the corresponding asymptotic
$\chi_{0.05}^2(3)$ distribution. For the sake of
brevity, we do not include the results for Mood's
and Klotz's normal tests, which are quite
similar. From Table 1, it can be seen that the
values are stable with respect to the choice of
sample sizes and parameters. We noted in our
Theorem 1 that the empirical rate of convergence
of the normalized random variable
$\Sigma_N^{-1/2} s_N$ to the $k$-variate normal
distribution $\mathcal{N}(0,I_k)$ depends on the
parameters of the GARCH process. The smaller the
parameters $\alpha_j$ and $\beta_j$, the faster
the convergence. This is intuitively clear
because larger values of $\alpha_j$ and $\beta_j$
imply not only more dependence, but also heavier
tails of the error distributions (cf. Basrak et
al. (2002)). More specifically, we observe that
the power of the tests for the DGP 1 is generally
higher than that
for the DGP 2 with respect to the sample sizes.\\
\indent Overall, the two bootstrap-based
statistics perform reasonably well in terms of
empirical size and power, and none of them
provides an obvious answer to the question of
what test statistic should be preferred.
Therefore, in practice we cannot know in advance
which of them would lead to a more powerful test.
Moreover, as the sample sizes and $\varphi$
increase, the size of both the tests converge to
the theoretical level and their powers generally
increase. When the error distributions are
sufficiently different, the power of the tests is
adequate for three different choices of the
sample size. It is worth noting that the highest
power of such tests is attained at $\varphi=1/3$.
  \begin{table}[htbp]
   \centering
  \caption{Proportion of rejections of $H_0$
for the bootstrap W and VdW tests at $r=5\%$}
\newcommand{\lw}[1]{\smash{\lower1.ex\hbox{#1}}}
\renewcommand{\arraystretch}{1.3}
\begin{tabular}{cccccccccc}
  \hline
 \multicolumn{10}{c}{ DGP 1}\\\hline
\lw{ $\varphi$}&&\multicolumn{2}{c}{
$n_1=n_2=n_3=100$}&& \multicolumn{2}{c}{
$n_1=n_2=n_3=300$}&&
\multicolumn{2}{c}{ $n_1=n_2=n_3=500$}\\
\cline{3-4} \cline{6-7} \cline{9-10}
&&\multicolumn{1}{c}{W}&\multicolumn{1}{c } {VdW
}&&\multicolumn{1}{c}{W} &\multicolumn{1}{c}{VdW}
&&\multicolumn{1}{c}{W } &\multicolumn{1}{c}{VdW}
\\
\hline 0&&0.045&0.046&&0.052&0.053 &&0.049&0.051
\\
1/9&&0.122&0.123&&0.171&0.169&&0.214&0.221
\\
1/5 &&0.321&0.331&&0.412&0.401 &&0.785&0.788
\\
 1/3&&0.821&0.818&&0.861&0.871&&0.913&0.912
\\
\hline
 \multicolumn{10}{c}{ DGP 2}\\\hline
0&&0.041&0.042&&0.056&0.055 &&0.045&0.046
\\
1/9&&0.102&0.104&&0.151&0.148&&0.193&0.195
\\
1/5 &&0.313&0.314&&0.393&0.401 &&0.712&0.717
\\
 1/3&&0.801&0.796&&0.815&0.817&&0.897&0.894
\\
\hline
 \end{tabular}
\end{table}

\indent We conclude this section with a simple
empirical example based on daily data. For this
purpose, we apply the bootstrap W and VdW tests
to the series of residuals obtained from the
estimation of a GARCH(1,1) on series of daily
individual stock returns for the three companies
(i) AMOCO, (ii) FORD and (iii) HP listed on New
York Stock Exchange. Each series starts from July
3, 1962, to December 31, 1991 with 7420
observations. In our analysis, however, we
consider the last 2000 data points from each
series from February 2, 1984, to December 31,
1991.\\
\indent Table 2 displays the empirical proportion
of rejections of $H_0$ for the W and VdW tests at
the 5\% significance level. The result shows that
the tests have similar desirable size and power
at the 5\% level. To this end, the results
provide enough evidence in support of the
simulation results. For all the three considered
series, the hypothesis of normality of the error
distributions is rejected at the 5\% level. The
bootstrap tests we studied in this paper have
reasonable size and can detect a misspecified
probability distribution of the errors in a GARCH
model with high probability.
    \begin{table}[htbp]
   \centering
  \caption{Proportion of rejections of $H_0$
for the bootstrap W and VdW tests at $r=5\%$}
\newcommand{\lw}[1]{\smash{\lower1.ex\hbox{#1}}}
\renewcommand{\arraystretch}{1.3}
\begin{tabular}{cccccccccc}
  \hline
$\varphi$&&0&&1/9&&1/5&&1/3\\
\hline
W&&0.050&&0.616&&0.981&&1.000\\
VdW&&0.049&&0.618&&0.978&&1.000\\\hline
\end{tabular}
 \end{table}

\section{Proof and Auxiliary Lemma}
In this section we provide Lemma 1 and the proof
of Theorem 1. Lemma 1 is useful for ordering
characteristic roots of a product of two
matrices (see e.g., Sen and Singer (1993)).\\
\\
\noindent\textbf{Lemma 1} (Courant). {\it Let $U$
and $V$ be positive semi-definite matrices.
Suppose that $V$ is nonsingular and that
$\underline{x}=(x_1,\ldots,x_k)^T\in
(-\infty,\infty)^k$ is a characteristic vector.
Then if the product $UV^{-1}$ is well defined,
and if $\upsilon_i$ denotes the $i$th
characteristic root of $UV^{-1}$ for
$i=1,\ldots,k$, we have
\[
\mbox{ch}_k(UV^{-1})=\upsilon_k=\inf_{\underline{x}}
\frac{\underline{x}^TU\underline{x}}{
\underline{x}^TV\underline{x}} \le
\sup_{\underline{x}}
\frac{\underline{x}^TU\underline{x}}
{\underline{x}^TV\underline{x}}=\upsilon_1=
\mbox{ch}_1(UV^{-1}).
\]  }
 \indent Next we provide the proof of Theorem
 1.\\

\noindent\textbf{Proof of Theorem 1.} Write
$d\hat{F}_{j,n_j}=d(\hat{F}_{j,n_j}-F_j+F_j)$ and
\begin{eqnarray*}
J\biggl(\frac{N}{N+1}\hat{\mathcal{H}}_N\biggl)&=&J(H_N)+
(\hat{\mathcal{H}}_N-H_N)J'(H_N)-\frac{\hat{\mathcal{H}}_N}
{N+1}J'(H_N)\nonumber\\
&&+\biggl[J\biggl(\frac{N}{N+1}
\hat{\mathcal{H}}_N\biggl)-J(H_N)-\biggl(\frac{N}{N+1}
\hat{\mathcal{H}}_N-H_N\biggl)J'(H_N)\biggl].
\end{eqnarray*}
Then the decomposition of (7) is given by
\[\hat{T}_{jN}=\mu_{jN}+B_{1N,j}
+B_{2N,j}+C_{1N,j}+C_{2N,j}+C_{3N,j},
\]
where
\begin{eqnarray*}
 B_{1N,j}&=&\int
J(H_N)d(\hat{F}_{j,n_j}-F_j)(x),\\
B_{2N,j}&=&\int
(\hat{\mathcal{H}}_N-H_N)J'(H_N)dF_j(x),\\
C_{1N,j}&=&-\frac{1}{N+1}\int
\hat{\mathcal{H}}_N J'(H_N)d\hat{F}_{j,n_j}(x),\\
C_{2N,j}&=&\int (\hat{\mathcal{H}}_N-H_N)J'(H_N)
d(\hat{F}_{j,n_j}-F_j)(x),\\
C_{3N,j}&=&\int
\biggl[J\biggl(\frac{N}{N+1}\hat{\mathcal{H}}_N
\biggl)-
J(H_N)\\
&&-\biggl(\frac{N}{N+1}\hat{\mathcal{H}}_N-H_N
\biggl)J'(H_N)\biggl] d\hat{F}_{j,n_j}(x).
\nonumber
\end{eqnarray*}
\indent To prove this theorem, it is necessary to
show that (i) the vector
$N^{1/2}(B_{1N,j}+B_{2N,j})_{1\le j\le k}$ when
properly normalized has a limiting Gaussian
distribution, and (ii) the $C_*$ terms are
uniformly of higher order. For (i), we observe
that the difference
$N^{1/2}(\hat{T}_{jN}-\mu_{jN}) _{1\le j\le k}
-N^{1/2}(B_{1N,j}+B_{2N,j})_{1\le j\le k}$ tends
to zero in probability and so the vectors
$N^{1/2}(\hat{T}_{jN}-\mu_{jN}) _{1\le j\le k}$
and $N^{1/2}(B_{1N,j}+B_{2N,j})_{1\le j\le k}$
possess the same limiting distribution.\\
\indent Let us now proceed to show the statement
(i). From (5), it is easily seen that
\begin{equation}
B_{1N,j} =\int J(H_N)d(F_{j,n_j}-F_j)(x)
+n_j^{-1/2}\mathcal{A}_j\int
J(H_N)d(xf_j(x))+o_p(1).
\end{equation}
Integrating $B_{2N,j}$ by parts, and using (6)
and (10), it follows that
\begin{eqnarray}
N^{1/2}(B_{1N,j}+B_{2N,j})
&=&N^{1/2}\biggl(-\sum_{{\scriptstyle i=1}
\atop{\scriptstyle i\ne j}}^k\lambda_{iN} \int
B_j(x)d(F_{i,n_i}-F_i)(x) \nonumber\\&& +\int
(J(H_N)-\lambda_jB_j(x))d(F_{j,n_j}-F_j)(x)
\nonumber\\
&&-n_j^{-1/2}\mathcal{A}_j\sum_{{\scriptstyle
i=1} \atop{\scriptstyle i\ne j}}^k\lambda_{iN}
\int x f_j(x)J'(H_N)dF_i(x) \nonumber\\&&
+\sum_{{\scriptstyle i=1} \atop{\scriptstyle i\ne
j}}^k\lambda_{iN}n_i^{-1/2}\mathcal{A}_i \int x
f_i(x)J'(H_N)dF_j(x)\biggl) +\, o_p(1)
\nonumber\\
&=& a_{jN}+b_{jN}+c_{jN}+d_{jN} + o_p(1)
\mbox{\qquad (say)},
\end{eqnarray}
where $ B_j(x)=\int_{x_0}^xJ'(H_N)dF_j(y) $ with $x_0$ determined
somewhat arbitrarily,
say by $H_N(x_0)=1/2$.\\
\indent In what follows, we shall first evaluate
the asymptotic variance of (11) and then the
asymptotic covariance to construct the dispersion
matrix $\Sigma_N$. For this purpose, first
consider $a_{jN}$ and write it as
 \begin{eqnarray*}
   \lefteqn{-N^{1/2}\lambda_{iN} \int
B_j(x)d(F_{i,n_i}-F_i)(x)}\\
& =&  N^{1/2}\lambda_{iN} \int (F_{i,n_i}-F_i)
J'(H_N)dF_j(x),\quad i=1,\ldots,
j-1,j+1,\ldots,k.
\end{eqnarray*}
 Then the mean is zero and the variance is
 \begin{eqnarray*}
  \lefteqn{E\biggl(N^{1/2}\lambda_{iN} \int
  (F_{i,n_i}-F_i)
J'(H_N)dF_j(x)\biggl)^2}\\
& =& E\biggl(N\lambda_{iN}^2 \int\!\!\!\! \int
(F_{i,n_i}(x)-F_i(x)) (F_{n_i}^{(i)}(y)-F_i(y))
J'(H_N(x))J'(H_N(y))dF_j(x)dF_j(y)\biggl)^2 \\
&=&  2\lambda_{iN}  \int\!\!\!\!\int
\limits_{\hspace*{-0.3cm} x<y}
F_i(x)(1-F_i(y))J'(H_N(x))J'(H_N(y))
dF_j(x)dF_j(y).
    \end{eqnarray*}
Note that the application of Fubini's theorem
permits the interchange of integral and expectation.\\
\indent By a similar argument, the variance of
 \[
b_{jN}=
 -N^{1/2}\sum_{{\scriptstyle i=1}
\atop{\scriptstyle i\ne j}}^k\lambda_{iN}
 \int  (F_{i,n_i}-F_i) J'(H_N)dF_j(x)
 \]
is given by
\begin{eqnarray*}
  \lefteqn{\frac{2}{\lambda_{jN}}  \sum_{{\scriptstyle i=1}
\atop{\scriptstyle i\ne j}}^k \lambda_{iN}^2 \int\!\!\!\!\int
\limits_{\hspace*{-0.3cm} x<y} F_j(x)(1-F_j(y))J'(H_N(x))J'(H_N(y))
dF_i(x)dF_i(y)} \\
&+& \frac{2}{\lambda_{jN}}  \sum_{{\scriptstyle
i,i'=1}\atop{\scriptstyle i\ne i',i\ne j,i'\ne j}}^k \lambda_{iN}
\lambda_{i'N} \int\!\!\!\!\int \limits_{\hspace*{-0.3cm} x<y}
F_j(x)(1-F_j(y))J'(H_N(x))J'(H_N(y))
dF_i(x)dF_{i'}(y)  \\
 &+& \frac{2}{\lambda_{jN}}  \sum_{{\scriptstyle
i,i'=1}\atop{\scriptstyle i\ne i',i\ne j,i'\ne
j}}^k \lambda_{iN} \lambda_{i'N} \int\!\!\!\!\int
\limits_{\hspace*{-0.3cm} y<x}
F_j(y)(1-F_j(x))J'(H_N(x))J'(H_N(y))
dF_i(x)dF_{i'}(y).
\end{eqnarray*}
Therefore, by observing that $a_{jN}$ and
$b_{jN}$ are mutually independent variables, it
follows by the result of Puri (1964) that
\begin{eqnarray}
\sigma_{1N,jj}&=&Var(a_{jN}+b_{jN})\nonumber\\
&=&2\Biggl\{\sum_{{\scriptstyle i=1}
\atop{\scriptstyle i\ne j}}^k \lambda_{iN}
\int\!\!\!\!\int \limits_{\hspace*{-0.3cm} x<y}
\Gamma_{iN}(x,y)dF_j(x)dF_j(y)\nonumber\\
&&+ \frac{1}{\lambda_{jN}}\sum_{{\scriptstyle
i=1} \atop{\scriptstyle i\ne j}}^k \lambda_{iN}^2
\int \!\!\!\!\int \limits_{\hspace*{-0.3cm}x<y}
\Gamma_{jN}(x,y)dF_i(x)dF_i(y)\Biggr\}\nonumber\\
&&+\frac{1}{\lambda_{jN}}\hspace*{-0.4cm}\sum_{{\scriptstyle
i,i'=1}\atop{\scriptstyle i\ne i',i\ne j,i'\ne
j}}^k \lambda_{iN}\lambda_{i'N} \Biggl\{\int
\!\!\!\!\int \limits_{\hspace*{-0.3cm}x<y}
\Gamma_{jN}(x,y)dF_i(x)dF_{i'}(y)
\nonumber\\
&& + \int \!\!\!\!\int
\limits_{\hspace*{-0.3cm}y<x}
\Gamma_{jN}(y,x)dF_i(x)dF_{i'}(y)\Biggr\},
\end{eqnarray}
where $\Gamma_{jN}(u,v)=F_j(u)(1-F_j(v))J'(H_N(u))J'(H_N(v))$. To
evaluate the same for $c_{jN}$ and $d_{jN}$, recall the result of
Francq and Zako\"ian (2004) that
\begin{eqnarray*}
Var(n_j^{1/2}(\hat{\theta}_{j,n_j}-\theta_{0j}))=
(\kappa_j-1)[\mathcal{U}(\theta_{0j})]^{-1},
\quad 1\le j\le k.
\end{eqnarray*}
In view of (5) and (11), it follows that
\begin{equation}
\sigma_{2N,jj} =Var(c_{jN})=
(\kappa_j-1)\omega_{jN}^T[\mathcal{U}(\theta_{0j})]^{-1}
\omega_{jN},
\end{equation}
 where $\omega_{jN}=-\lambda_{jN}^{-1/2}
 \sum\limits_{{\scriptstyle
i=1} \atop{\scriptstyle i\ne j}}^k \lambda_{iN} \int x f_j(x)J'(H_N)
dF_i(x)\times \tau_j$ with
$\tau_j=(\tau_{j,1},\ldots,\tau_{j,p_j+q_j+1})^T$, and analogously
 \begin{equation}
\sigma_{3N,jj} =Var(d_{jN})= \sum_{{\scriptstyle
i=1} \atop{\scriptstyle i\ne j}}^k
(\kappa_i-1)\nu_{iN}^T[\mathcal{U}(\theta_{0i})]^{-1}
\nu_{iN},
\end{equation}
where $\nu_{iN}=\lambda_{iN}^{1/2} \int x f_i(x)
J'(H_N)dF_j(x)\times \tau_i$. Moreover, by
independence of $X_{j,1},\ldots,X_{j,n_j}$, it
remains to evaluate
\[
K_{1N,j}=2E(a_{jN}d_{jN})\quad \mbox{and}\quad
K_{2N,j}=2E(b_{jN}c_{jN}).
\]
Using (11), we obtain
\[
K_{1N,j}=2\sum_{{\scriptstyle i=1}
\atop{\scriptstyle i\ne
j}}^k\lambda_{iN}\int\!\!\!\!\int
E[(n_i^{1/2}(F_{i,n_i}(x)-F_i(x))\mathcal{A}_i]
\psi_{iN}(x,y)dF_j(x)dF_j(y),
\]
where
$\psi_{iN}(u,v)=vf_i(v)J'(H_N(u))J'(H_N(v))$. To
obtain an explicit expression of $K_{1N,j}$, it
is necessary to evaluate $E[\cdot]$. From the
result of Berkes and Horv\'ath (2003) and (5), we
find that
\[
E[n_i^{1/2}(F_{i,n_i}(x)-F_i(x))\mathcal{A}_i]
=\sum_{l=1}^{p_i+q_i+1}\tau_{i,l}h_i^l(x),
\]
where $h_i^l(v)=\delta_i^l\int_{-\infty}^v
(u^2-1)f_i(u)du$ with $\delta _i^{l}=E
(Z_{t}^l(\theta_{0i}))$, $1\le i\le k$. Then
\[
K_{1N,j}=2\sum_{{\scriptstyle i=1}
\atop{\scriptstyle i\ne
j}}^k\sum_{l=1}^{p_i+q_i+1}\lambda_{iN}\tau_{i,l}
\int\!\!\!\!\int
h_i^l(x)\psi_{iN}(x,y)dF_j(x)dF_j(y),
\]
and similarly
\[
K_{2N,j}=\frac{2}{\lambda_{jN}}\sum_{{\scriptstyle
i=1} \atop{\scriptstyle i\ne
j}}^k\sum_{l=1}^{p_j+q_j+1}\lambda_{iN}^2\tau_{j,l}\int\!\!\!\!\int
h_j^l(x)\psi_{jN}(x,y)dF_i(x)dF_i(y).
\]
Therefore, the variance terms when combined yield
\begin{equation}
   \sigma_{N,jj}=\sigma_{1N,jj}+\sigma_{2N,jj}
+\sigma_{3N,jj}+\gamma_{N,jj},
\end{equation}
where $\gamma_{N,jj}= K_{1N,j}+K_{2N,j}$.\\
\indent We next turn to evaluate the covariance
terms. For this purpose, rewrite (11) as
\begin{eqnarray*}
\lefteqn{ N^{1/2}(B_{1N,j}+B_{2N,j})
}\nonumber\\
&=& N^{1/2}\sum_{i=1}^k\lambda_{iN} \biggl( -\int
(F_{j,n_j}(x)- F_j(x))J'(H_N)dF_i(x)
\nonumber\\&& +\int
(F_{n_i}^{(i)}(x)-F_i(x))J'(H_N)dF_j(x)
-n_j^{-1/2}\mathcal{A}_j\int
 x f_j(x)J'(H_N)dF_i(x)\nonumber\\
&&+n_i^{-1/2} \mathcal{A}_i\int x
f_i(x)J'(H_N)dF_j(x)\biggl)
+ \,o_p(1)\nonumber\\
&=& a_{1N,j}+b_{1N,j}+c_{1N,j}+d_{1N,j} + o_p(1),
\quad\mbox{ (say)}.
\end{eqnarray*}
By independence of $X_{j,1},\ldots,X_{j,n_j}$,
$1\le j\le k$, we first compute
\begin{eqnarray*}
\sigma_{1Njj'}&=&Cov(a_{1N,j}+b_{1N,j},a_{1N,j'}
+b_{2N,j'})\\
&=&E(a_{1N,j}b_{1N,j'})+E(b_{1N,j}a_{1N,j'})
+E(b_{1N,j}b_{1N,j'}),\quad j\not=j'=1,\ldots, k.
\end{eqnarray*}
From
\begin{eqnarray*}
 a_{1N,j}b_{1N,j'} &=&-N \sum_{i=1}^k\sum_{l=1}^k
 \lambda_{iN}\lambda_{lN}\int \!\!\!\!\int
(F_{j,n_j}(x)-F_j(x))(F_{l,n_l}(y)-F_l(y))\\
&&\times J'(H_N(x))J'(H_N(y))dF_i(x)dF_{j'}(y),
\end{eqnarray*}
it follows by using again the result of Puri
(1964) that
 \begin{eqnarray*}
 E(a_{1N,j}b_{1N,j'}) &=&- \sum_{i=1}^k
 \lambda_{iN}\int \!\!\!\!\int
\limits_{\hspace*{-0.3cm}x<y} F_j(x)(1-F_j(y))
 J'(H_N(x))J'(H_N(y))dF_i(x)dF_{j'}(y)   \\
 &&- \sum_{i=1}^k
 \lambda_{iN}\int \!\!\!\!\int
\limits_{\hspace*{-0.3cm}y<x} F_j(y)(1-F_j(x))
 J'(H_N(x))J'(H_N(y))dF_i(x)dF_{j'}(y).
\end{eqnarray*}
 In the same way, we have
  \begin{eqnarray*}
 E(b_{1N,j}a_{1N,j'}) &=&- \sum_{i=1}^k
 \lambda_{iN}\int \!\!\!\!\int
\limits_{\hspace*{-0.3cm}x<y}
F_{j'}(x)(1-F_{j'}(y))
 J'(H_N(x))J'(H_N(y))dF_i(x)dF_{j}(y)   \\
 &&- \sum_{i=1}^k
 \lambda_{iN}\int \!\!\!\!\int
\limits_{\hspace*{-0.3cm}y<x}
F_{j'}(y)(1-F_{j'}(x))
 J'(H_N(x))J'(H_N(y))dF_i(x)dF_{j}(y)
\end{eqnarray*}
and
 \begin{eqnarray*}
 E(b_{1N,j}b_{1N,j'}) &=& \sum_{i=1}^k
 \lambda_{iN}\int \!\!\!\!\int
\limits_{\hspace*{-0.3cm}x<y} F_i(x)(1-F_i(y))
 J'(H_N(x))J'(H_N(y))dF_j(x)dF_{j'}(y)   \\
 &&+ \sum_{i=1}^k
 \lambda_{iN}\int \!\!\!\!\int
\limits_{\hspace*{-0.3cm}y<x}
F_{i}(y)(1-F_{i}(x))
 J'(H_N(x))J'(H_N(y))dF_j(x)dF_{j'}(y).
\end{eqnarray*}
Therefore,
    \begin{eqnarray*}
    \sigma_{1Njj'} &=& -\sum_{i=1}^k\lambda_{iN}
    \biggl(\int \!\!\!\!\int
\limits_{\hspace*{-0.3cm}x<y}
\Gamma_{jN}(x,y)dF_i(x)dF_{j'}(y)+\int
\!\!\!\!\int \limits_{\hspace*{-0.3cm}y<x}
\Gamma_{jN}(y,x)dF_i(x)dF_{j'}(y)\biggl)\nonumber\\
&&-\sum_{i=1}^k\lambda_{iN}\biggl(\int
\!\!\!\!\int \limits_{\hspace*{-0.3cm}x<y}
\Gamma_{j'N}(x,y)dF_i(x)dF_{j}(y)+\int
\!\!\!\!\int \limits_{\hspace*{-0.3cm}y<x}
\Gamma_{j'N}(y,x)dF_i(x)dF_{j}(y)
\biggl)\nonumber\\
&&+\sum_{i=1}^k\lambda_{iN}\biggl(\int
\!\!\!\!\int \limits_{\hspace*{-0.3cm}x<y}
\Gamma_{iN}(x,y)dF_j(x)dF_{j'}(y)+\int
\!\!\!\!\int \limits_{\hspace*{-0.3cm}y<x}
\Gamma_{iN}(y,x)dF_j(x)dF_{j'}(y)\biggl).
 \end{eqnarray*}
Now we turn to evaluate, for $j\not=j'$,
\begin{eqnarray*}
L_{1N,jj'}
=E(a_{1N,j}d_{1N,j'})+E(d_{1N,j}a_{1N,j'})\quad\mbox{and}\quad
L_{2N,jj'}=E(b_{1N,j}c_{1N,j'})+E(c_{1N,j}b_{1N,j'}).
\end{eqnarray*}
In analogy with the preceding $K_{*}$ terms, we
have
\begin{eqnarray*}
L_{1N,jj'}&=&-\sum_{i=1}^k\lambda_{iN}\biggl(\sum_{l=1}^
{p_j+q_j+1}\tau_{j,l}\int\!\!\!\!\int
h_j^l(x)\psi_{j}(x,y)dF_i(x)dF_{j'}(y)\nonumber\\
&&+\sum_{l=0}^{p_{j'}+q_{j'}+1}\tau_{j',l}\int\!\!\!\!\int
h_{j'}^l(x)\psi_{j'}(x,y)dF_i(x)dF_{j}(y)\biggl)\\
\end{eqnarray*}
and
\begin{eqnarray*}
L_{2N,jj'}&=&-\sum_{i=1}^k\lambda_{iN}\biggl(
\sum_{l=1}^{p_j+q_j+1}\tau_{j,l}\int\!\!\!\!\int
h_{j}^l(y)\psi_{j}(y,x)dF_i(x)dF_{j'}(y)\nonumber\\
&&+
\sum_{l=1}^{p_{j'}+q_{j'}+1}\tau_{j',l}\int\!\!\!\!\int
h_{j'}^l(y)\psi_{j'}(y,x)dF_i(x)dF_{j}(y)
\biggl).
\end{eqnarray*}
Therefore, combining the covariance terms
produces
\begin{equation}
   \sigma_{N,jj'}=  \sigma_{1N,jj'}+  \sigma_{2N,jj'},
\end{equation}
where $\sigma_{2N,jj'}=L_{1N,jj'}+  L_{2N,jj'}$, $j\ne j'=1,\ldots,k$.\\
\indent Hence, using (13)$-$(16) and the central
limit theorems for martingale differences given
by Berkes and Horv\'ath (2003), and Francq and
Zako\"ian (2004), we may conclude that
\begin{eqnarray*}
N^{1/2}\Sigma_N^{-1/2}(B_{1N,j}+B_{2N,j})_{1\le
j\le k}\stackrel {d}{\longrightarrow}{\mathcal
N}(0,I_k) \quad \mbox{as $N\to\infty$}.
\end{eqnarray*}

\indent Next, we turn to show statement (ii). For
this purpose, we require the following elementary
results (see Puri (1964)).
\begin{itemize}
 \item[(i)] $H_N\ge \lambda_{jN}F_j\ge \lambda_0F_j$,
 \quad $1\le j\le k$.
 \item[(ii)] $(1-F_j)\le (1-H_N)/\lambda_{jN}\le
 (1-H_N)/\lambda_0$,\quad $1\le j\le k$.
 \item[(iii)] $F_j(1-F_j)\le H_N(1-H_N)/
 \lambda_{jN}^2
 \le H_N(1-H_N)/\lambda_0^2$,\quad $1\le j\le k$.
 \item[(iv)] $dH_N\ge \lambda_{jN}dF_j\ge\lambda_0
 dF_j$,
 \quad $1\le j\le k$.
\item[(v)] Let ($\vartheta_{1N},\vartheta_{2N}$)
be the interval $S_{N_\epsilon}$ such that
\begin{equation}
S_{N_\epsilon}=\{x:H_N(1-H_N)>\eta_\epsilon \lambda_0/N\},
\end{equation}
\end{itemize}
where $\epsilon>0$ is arbitrarily small and
$\eta_\epsilon(>0)$ depends $\epsilon$. Thus,
\begin{equation}
  \eta_\epsilon< N (H_N(\vartheta_{1N})),\qquad
 (1-H_N(\vartheta_{2N}))<\eta_\epsilon(1+N^{-1}\eta_\epsilon).
\end{equation}
Hence, $\eta_\epsilon$ can be chosen
independently of $F_j$ and $\lambda_{jN}$ in such
a way that
\begin{equation}
   N(H_N(\vartheta_{1N})+(1-H_N(\vartheta_{2N})))\le
\epsilon.
\end{equation}
From (19), it follows that
\begin{eqnarray}
  \lefteqn{
  P(\varepsilon_{j,t}\in S_{N_\epsilon},1\le t\le
n_j, 1\le j\le k)}
\nonumber\\
 &&= P(\vartheta_{1N}\le
\varepsilon_{j,t}\le \vartheta_{2N})\nonumber \\
&&=\prod_{j=1}^N[H_j(\vartheta_{2N})-
H_j(\vartheta_{1N})]\nonumber \\
&&= \prod_{j=1}^N\{1-[H_j(\vartheta_{1N})+1-
H_j(\vartheta_{2N})]\}\nonumber\\
&&\ge 1-\sum_{j=1}^N [H_j(\vartheta_{1N})+1-
H_j(\vartheta_{2N})]\nonumber \\
&&=1-N[H_N(\vartheta_{1N})+(1-
H_N(\vartheta_{2N}))]\nonumber\\
&&\ge 1- \epsilon.
\end{eqnarray}

\indent Let us first evaluate $C_{1N,j}$. By (6)
and $d\hat{F}_{j,n_j}=d(\hat{F}_{j,n_j}-F_j
+F_j)$, we have
\begin{eqnarray}
C_{1N,j} &=&\frac{-1}{N+1}\int
\mathcal{H}_NJ'(H_N)dF_{j,n_j}(x)\nonumber\\
&&- \frac{1}{N(N+1)}\sum_{i=1}^kn_i^{1/2}\mathcal{A}_i\int
xf_i(x)J'(H_N)dF_{j,n_j}(x)\nonumber\\
&&-\frac{n_j^{-1/2}}{N+1}\mathcal{A}_j\int
\mathcal{H}_NJ'(H_N)d(xf_j(x))\nonumber\\
&&-\frac{n_j^{-1/2}\mathcal{A}_j}
{N(N+1)}\sum_{i=1}^kn_i^{1/2}\mathcal{A}_i\int
xf_i(x)J'(H_N)d(xf_j(x))
+\,o_p(N^{-1})\nonumber\\
&=&\sum_{i=1}^4 C_{1iN,j}
+o_p(N^{-1})\nonumber,\quad\mbox{(say)}.
\end{eqnarray}
The proof of $C_{11N,j}=o_p(N^{-1/2})$ follows
precisely the same arguments as in Puri (1964).
Next we turn to $C_{12N,j}$. By (A.2) and (A.3),
we obtain
\begin{eqnarray*}
   |C_{12N,j}|&\le&
\frac{1}{N}\sum_{i=1}^kn_i^{1/2}|
\mathcal{A}_i|\frac{1}{N+1}\int
|xf_i(x)||J'(H_N)|dF_{j,n_j}(x)\\
&\le& \frac{1}{N}\sum_{i=1}^kn_i^{1/2}|
\mathcal{A}_i|\frac{1}{N}\int |J'(H_N)|
dF_{j,n_j}(x).
\end{eqnarray*}
In a similar fashion as the proof for
$C_{11N,j}$, it follows that
\[
\frac{1}{N}\int |J'(H_N)|
dF_{j,n_j}(x)=o_p(N^{-1/2}),
\]
which, combined with the fact
\begin{equation}
\frac1N\sum_{i=1}^k
n_i^{1/2}|\mathcal{A}_i|=\mathcal{O}_p
\biggl(\frac1N\sum_{i=1}^kn_i^{1/2}\biggl),
\end{equation}
implies $ C_{12N,j}=o_p(N^{-1})$. Next consider
\[
  C_{13N,j}=-n_j^{-1/2}\mathcal{A}_j(C_{13N,j}^*
  +C_{13N,j}^{**}),
\]
 where
 \begin{eqnarray*}
   C_{13N,j}^*  &=& \frac{1}{N+1}\int _{S_{N_\epsilon}}
   \mathcal{H}_N J'(H_N) d(xf_j(x)),    \\
   C_{13N,j}^{**} &=& \frac{1}{N+1}\int _{S_{N_\epsilon}^c}
   \mathcal{H}_N J'(H_N) d(xf_j(x))
 \end{eqnarray*}
and $S_{N_\epsilon}^c$ is the complementary event
of $S_{N_\epsilon}$. Let us first deal with
$C_{13N,j}^*$. In view of (A.2), (A.3), (17) and
(18), it follows that
\begin{eqnarray}
| C_{13N,j}^* | &\le& \frac{K}{c_jN} \int _{S_{N_\epsilon}}
|J'(H_N)|
dF_j(x)\nonumber\\
&\le& \frac{K}{c_jN} \int _{S_{N_\epsilon}}
   [H_N(1-H_N)]^{-\frac32+\delta}dH_N(x)\nonumber\\
   &\le & \frac{K}{c_jN} \int_{\frac{K}{N}}^1
   H_N^{-\frac32+\delta}dH_N(x) \nonumber\\
   &\le & \frac{K}{N^{\frac12+\delta}}.
\end{eqnarray}
Now using the Markoff inequality, we obtain
\[
P(| C_{13N,j}^* |>mN^{-1/2})\le
\frac{K}{N^{\frac12+\delta}}
\frac{N^{1/2}}{m}=\frac{K}{mN^\delta},
\]
where $m>0$ and $K$ may depend on $\epsilon$.
Next consider $C_{13N,j}^{**}$. Write
$H_1=H_N(\vartheta_{1N})$ and
$H_2=H_N(\vartheta_{2N})$. Then from (17) and
(18), we have $H_1=1-H_2<K/N$. By (20), we are
certain that $\varepsilon_{j,t}\not\in
S_{N_\epsilon}^c$ and
 \begin{eqnarray}
| C_{13N,j}^{**} | &\le& \frac{K}{c_jN}\biggl(
\int _0^{H_1} [H_N(1-H_N)]^{-\frac32+\delta}
dH_N(x)\nonumber\\
&&+ \int _{H_2}^1 [H_N(1-H_N)]^{-\frac32+\delta}
dH_N(x)\biggl)
\nonumber\\
&\le& \frac{K}{c_jN} \int _0^{H_1}
H_N^{-\frac32+\delta} dH_N(x)\nonumber\\
   &\le & \frac{K}{N^{\frac12+\delta}}.
\end{eqnarray}
Therefore, by using (21), we have
\begin{equation}
  C_{13N,j} =o_p(N^{-1/2}).
\end{equation}
Similarly, it can be shown that
$C_{14N,j}=o_p(N^{-1})$. Consequently, we have
\[C_{1N,j}=o_p(N^{-1/2}).
\]
\indent Next, we consider $C_{2N,j}$. By analogy
with the first $C$ term, we have
\begin{eqnarray*}
C_{2N,j} &=&
\int(\mathcal{H}_N-H_N)J'(H_N)d(F_{j,n_j}-F_j)(x)
\\&&
+\frac{1}{N}\sum_{i=1}^kn_i^{1/2}\mathcal{A}_i
\int xf_i(x)J'(H_N)d(F_{j,n_j}-F_j)(x) \\
&&+\frac{n_j^{-1/2}\mathcal{A}_j}{N}
\sum_{i=1}^kn_i^{1/2}\mathcal{A}_i\int xf_i(x)J'(H_N)d(xf_j(x))
\\&&
+n_j^{-1/2}\mathcal{A}_j\int (\mathcal{H}_N-H_N)J'(H_N)d(xf_j(x))
+o_p(N^{-1})\\
&=&\sum_{i=1}^4C_{2iN,j}
+o_p(N^{-1}),\quad\mbox{(say)}.
\end{eqnarray*}
The proof of $C_{21N,j}=o_p(N^{-1/2})$ is
identical to that of Puri (1964). Next, we
consider
 \[
  C_{22N,j}=\frac1N\sum_{i=1}^k n_j^{1/2}\mathcal{A}_i
  (C_{22N,j}^*
  +C_{22N,j}^{**}),
\]
for which, it suffices to show
\begin{eqnarray}
  C_{22N,j}^* &=& \int _{S_{N_\epsilon}}
xf_i(x)J'(H_N)d(F_{j,n_j}-F_j)(x)=o_p(1), \\
  C_{22N,j}^{**} &=& \int _{S_{N_\epsilon}^c}
xf_i(x)J'(H_N)d(F_{j,n_j}-F_j)(x)=o_p(1).
\end{eqnarray}
Note that from (A.2) and (A.3), we can find $K>0$
such that $ |xf_j(x)|\le KH_N(1-H_N)$. Then from
(17), (18) and (22), it follows that (25) is
dominated by
 \begin{eqnarray*}
| C_{22N,j}^* | &\le& K \int
_{S_{N_\epsilon}}|xf_j(x)| |J'(H_N)|
|d(F_{j,n_j}-F_j)(x)|\nonumber\\
&\le& K \int _{S_{N_\epsilon}}
   [H_N(1-H_N)]^{-\frac12+\delta}|d(F_{j,n_j}-F_j)(x)|\nonumber\\
   &\le & n_j^{-1/2} \int_{\frac{K}{N}}^1
   \mathcal{O}(N^{\frac12-\delta})|d[n_j^{1/2}(F_{j,n_j}
   -F_j)(x)]|
   =o_p(1).
\end{eqnarray*}
Likewise, it is easy to show from (23) that (26)
is dominated by
  \begin{eqnarray*}
| C_{22N,j}^{**} | &\le& K \biggl(\int _0^{H_1}
[H_N(1-H_N)]^{-\frac12+\delta}
|d(F_{j,n_j}-F_j)(x)|\\
&&+
    \int _{H_2}^1
[H_N(1-H_N)]^{-\frac12+\delta}
|d(F_{j,n_j}-F_j)(x)|
\biggl)\nonumber\\
   &\le & n_j^{-1/2} \int_0^{H_1}
   \mathcal{O}(N^{\frac12-\delta})|d[n_j^{1/2}(F_{j,n_j}
   -F_j)(x)]|
   =o_p(1).
\end{eqnarray*}
Therefore, it follows from (21) that
$C_{22N,j}=o_p(N^{-1/2})$. The proof for
$C_{23N,j}=o_p(N^{-1/2})$ is analogous to (24).
To complete the assertion for $C_{2N,j}$, it
remains to evaluate $C_{24N,j}= \mathcal{A}_j
  (C_{24N,j}^*
  +C_{24N,j}^{**})$,
where
 \begin{eqnarray*}
   C_{24N,j}^* &=& n_j^{-1/2}\int
_{S_{N_\epsilon}}
(\mathcal{H}_N-H_N)J'(H_N)d(xf_j(x)), \\
   C_{24N,j}^{**} &=& n_j^{-1/2}\int
_{S_{N_\epsilon}^c} (\mathcal{H}_N-H_N)J'(H_N)d(xf_j(x)).
 \end{eqnarray*}
By virtue of Puri and Sen (1993, Theorem 2.11.10), write
\begin{eqnarray}
{\cal I}_N(\delta')=\sup_x  \frac{  N^{1/2}
|\mathcal{H}_N(x)-H_N(x)|} {[H_N(x)(1-H_N(x))]^{\frac12-\delta'}}\le
C^*,\quad \delta'>0,\quad C^*>0
\end{eqnarray}
so that $P({\cal I}_N(\delta'))\ge 1-\epsilon$.
Then, if we let $\delta'<\delta$, it follows from
(A.2)$-$(A.4), (22) and (27) that
\begin{eqnarray*}
C_{24N,j}^* &=& \frac{n_j^{-1/2}}{c_j}\int _{S_{N_\epsilon}}
|\mathcal{H}_N-H_N||J'(H_N)|dF_j(x)\\
&\le& \frac{Kn_j^{-1/2}}{c_j}\int
_{S_{N_\epsilon}} \mathcal{O}(N^{-1/2})
[H_N(1-H_N)]^{\delta-\delta'-1} dH_N(x)\nonumber\\
&\le&
\frac{Kn_j^{-1/2}}{c_j}\mathcal{O}(N^{-1/2})
  \int_{\frac{K}{N}}^1
  H_N^{\delta-\delta'-1}dH_N(x)\\
   &=&
   \mathcal{O}(N^{\delta'-\delta-1})=o(N^{-1})
\end{eqnarray*}
and similarly from (23) that
 \begin{eqnarray*}
C_{24N,j}^{**} &\le&
\frac{Kn_j^{-1/2}}{c_j}\mathcal{O}(N^{-1/2})
  \int_0^{H_1}
  H_N^{\delta-\delta'-1}dH_N(x)
   = o(N^{-1}).
\end{eqnarray*}
 Hence, $C_{24N,j}=o_p(N^{-1/2})$. Consequently,
 we have
$$C_{2N,j}=o_p(N^{-1/2}).$$
\indent Finally, we evaluate $C_{3N,j}$.
Following the preceding $C_*$ term, and using
\begin{eqnarray*}
 J\biggl(\frac{N}{N+1}\hat{\mathcal{H}}_N\biggl)
&=&J(H_N)+\biggl(\frac{N}{N+1}\hat{\mathcal{H}}_N-H_N\biggl)\\
&&\times J'\biggl(\varrho
H_N+(1-\varrho)\frac{N}{N+1}\hat{\mathcal{H}}_N\biggl),\quad
0<\varrho<1,
\end{eqnarray*}
we obtain
\begin{eqnarray*}
C_{3N,j}
&=&\int\biggl(\frac{N}{N+1}\mathcal{H}_N-H_N\biggl)
\\&&\times
\biggl[J'\biggl(\varrho
H_N+(1-\varrho)\frac{N}{N+1}\hat{\mathcal{H}}_N\biggl)
-J'(H_N)\biggl]dF_{j,n_j}(x)
\\&&
+\frac{1}{N+1}\sum_{i=1}^k
n_i^{1/2}\mathcal{A}_i\int x f_i(x)
\\&&\times
\biggl[J'\biggl(\varrho H_N+(1-\varrho)
\frac{N}{N+1}\hat{\mathcal{H}}_N\biggl)-J'(H_N)
\biggl]dF_{j,n_j}(x)\\
&&+n_j^{-1/2}\mathcal{A}_j
\int\biggl(\frac{N}{N+1}\mathcal{H}_N-H_N\biggl)
\\&&\times
\biggl[J'\biggl(\varrho H_N+(1-\varrho)
\frac{N}{N+1}\hat{\mathcal{H}}_N\biggl)
-J'(H_N)\biggl]d(xf_j(x))\\
&&+\frac{n_j^{-1/2}\mathcal{A}_j}{N+1}
\sum_{i=1}^kn_i^{1/2}\mathcal{A}_i\int x f_i(x)\\
&&\times\biggl[J'\biggl(\varrho
H_N+(1-\varrho)\frac{N}{N+1}\hat{\mathcal{H}}_N
\biggl)-J'(H_N)\biggl]d(xf_j(x))
+o_p(N^{-1})\\
&=&\sum_{i=1}^4C_{3iN,j}
+o_p(N^{-1}),\quad\mbox{(say)}.
\end{eqnarray*}
First consider
$C_{31N,j}=C_{31N,j}^*+C_{31N,j}^{**}$, where
\begin{eqnarray*}
  C_{31N,j}^* &=&
  \int_{S_{N_\epsilon}}\biggl(\frac{N}{N+1}
  \mathcal{H}_N-H_N\biggl)
\\
&&\times \biggl[J'\biggl(\varrho
H_N+(1-\varrho)\frac{N}{N+1}\hat{\mathcal{H}}_N\biggl)
-J'(H_N)\biggl]dF_{j,n_j}(x),\\
  C_{31N,j}^{**} &=& \int_{S_{N_\epsilon}^c}\biggl(\frac{N}{N+1}\mathcal{H}_N-H_N\biggl)
\\&&\times \biggl[J'\biggl(\varrho
H_N+(1-\varrho)\frac{N}{N+1}\hat{\mathcal{H}}_N\biggl)
-J'(H_N)\biggl]dF_{j,n_j}(x).
\end{eqnarray*}
 To evaluate $C_{31N,j}^*$, first
 note from (6), (A.2), (A.3) and (21) that
\begin{eqnarray}
  \lefteqn{H_N-\biggl(\varrho
H_N+(1-\varrho)\frac{N}{N+1}
\hat{\mathcal{H}}_N\biggl)}\nonumber\\
&=&(1-\varrho )\biggl(H_N-\frac{N}{N+1}
\hat{\mathcal{H}}_N\biggl)  \nonumber\\
 &=& (1-\varrho) \biggl[\biggl(H_N-\frac{N}{N+1}
 \mathcal{H}_N\biggl)
 -\frac{N}{N+1}\sum_{j=1}^kn_j^{-1/2}
 \lambda_{jN}\mathcal{A}_jxf_j(x)\biggr]
 +o_p(N^{-1/2})\nonumber\\
   &=& N^{-1/2}\mathcal{O}(1)
   \{ 1+[H_N(1-H_N)]^{\frac12-\delta'}\},
\end{eqnarray}
where $\mathcal{O}(1)$ is uniform in $x$. Then
from (18) and (28), it follows that
   \begin{eqnarray*}
     \lefteqn{1-\biggl(\varrho
H_N+(1-\varrho)\frac{N}{N+1}
\hat{\mathcal{H}}_N\biggl)}\\
&=& (1-H_N)+ (1-\varrho )\biggl(H_N-\frac{N}{N+1}
\hat{\mathcal{H}}_N\biggl) \\
&=&(1-H_N) + \mathcal{O}(N^{-1/2})\{
1+[H_N(1-H_N)]^{\frac12-\delta'}\}\\
&=&(1-H_N)\{ 1+ \mathcal{O}(N^{-1/2})
H_N^{\frac12-\delta'}(1-H_N)^{-\frac12-\delta'}\}+
\mathcal{O}(N^{-1/2})\\
 &=&(1-H_N)\{ 1+ \mathcal{O}(N^{-1/2})
\mathcal{O}(N^{\frac12+\delta')})\}
+\mathcal{O}(N^{-1/2})\\
&=&(1-H_N)( 1+ o(1))
+\mathcal{O}(N^{-1/2})\\
&=& (1-H_N)+ \mathcal{O}(N^{-1/2})
\end{eqnarray*}
or equivalently
\[
\biggl(1-\biggl(\varrho
H_N+(1-\varrho)\frac{N}{N+1}\hat{\mathcal{H}}_N
\biggl)\biggl) \times (1-H_N)^{-1}=
1+\mathcal{O}(N^{-1/2}).
\]
Thus, for sufficiently large $N>0$, we can find
$\zeta>0$ such that
\begin{equation}
\inf_x\left(\varrho
H_N(x)+(1-\varrho)\frac{N}{N+1}\hat{\mathcal{H}}_N(x)\right)
\left(\frac{1-(\varrho H_N(x)+(1-\varrho)
\frac{N}{N+1}\hat{\mathcal{H}}_N(x))}{H_N(x)(1-H_N(x))}\right)
>\zeta
\end{equation}
with probability $\ge 1-\epsilon$. Now write
\begin{eqnarray}
|C_{31N,j}^*|&\le& \int \biggl|
\frac{N}{N+1}\mathcal{H}_N-H_N\biggl|
\nonumber\\
&&\times \biggl|J'\biggl[\varrho H_N+
(1-\varrho)\frac{N}{N+1}\hat{\mathcal{H}}_N\biggl]-
J'(H_N)\biggl|dF_{j,n_j}(x)\nonumber\\
&=& \int{\cal Q}_NdF_{j,n_j}(x),\quad \mbox{
(say)}.
\end{eqnarray}
Then it is easy to show from (A.2), (22), (27),
(29) and (30) that
\begin{eqnarray}
E\int _{S_{N_\epsilon}}{\cal
Q}_NdF_{j,n_j}(x)&\le&
K(1+\zeta^{\delta-\frac32})\mathcal{O}(N^{-1/2})
\nonumber\\
&&\times\int_{S_{N_\epsilon}}
[H_N(1-H_N)]^{\delta-\delta'-1}dH_N(x) \nonumber\\
&\le
&K(1+\zeta^{\delta-\frac32})\mathcal{O}(N^{-1/2})
 \int_{\frac{K}{N}}^1 H_N^{\delta-\delta'-1}
dH_N(x).
\end{eqnarray}
Thus, $\mathcal{Q}_N(x)$ is integrable and
converges to 0 in probability. Hence, by virtue
of the dominated convergence theorem and (31), it
is seen that $C_{31N,j}^*=o_p(N^{-1/2})$.
Similarly, we can show
$C_{31N,j}^{**}=o_p(N^{-1/2})$ by using the
arguments of (23) and (31). Next consider
\[C_{32N,j}=\frac{1}{N+1}\sum_{i=1}^k
n_i^{1/2}\mathcal{A}_i(C_{32N,j}^*+C_{32N,j}^{**}),
\]
  where
\begin{eqnarray*}
  C_{32N,j}^* &=& \int_{S_{N_\epsilon}}  x f_i(x)
\biggl\{J'\biggl[\varrho H_N+(1-\varrho)
\frac{N}{N+1}\hat{\mathcal{H}}_N\biggl]-J'(H_N)
\biggl\}dF_{j,n_j}(x) \\
  C_{32N,j}^{**} &=& \int_{S_{N_\epsilon}^c}  x f_i(x)
\biggl\{J'\biggl[\varrho H_N+(1-\varrho)
\frac{N}{N+1}\hat{\mathcal{H}}_N\biggl]-J'(H_N)
\biggl\}dF_{j,n_j}(x).
\end{eqnarray*}
Let us first evaluate $ C_{32N,j}^* $. Recalling
$|xf_j(x)|\le KH_N(1-H_N)$, and using the
arguments of $C_{31N,j}^*$ and (A.2), we obtain
\begin{eqnarray}
E(|C_{32N,j}^*|) &\le& \int_
{S_{N_\epsilon}}|xf_i(x)| \biggl|J'\biggl[\varrho
H_N+(1-\varrho)\frac
{N}{N+1}\hat{\mathcal{H}}_N\biggl]-J'(H_N)
\biggl|dF_j(x)\nonumber\\
&\le& K(1+\zeta^{\delta-\frac32})
\int_{S_{N_\epsilon}}
[H_N(1-H_N)]^{\delta-\frac12}dH_N(x)\nonumber\\
  &\le & K(1+\zeta^{\delta-\frac32})
  \int_{\frac{K}{N}}^1
H_N^{\delta-\frac12}dH_N(x).
\end{eqnarray}
In analogy with (23) and (32), we can show
$C_{32N,j}^{**}=o_p(N^{-1/2})$. Hence, from (21),
we have $C_{32N,j}=o_p(N^{-1/2})$. Next, we
evaluate $C_{33N,j}=C_{33N,j}^* +C_{33N,j}^{**}$,
where
\begin{eqnarray*}
  C_{33N,j}^* &=& n_j^{-1/2}\mathcal{A}_j\int_{S_{N_\epsilon}}\biggl(\frac{N}{N+1}\mathcal{H}_N-H_N\biggl)
\\&&\times
\biggl\{J'\biggl[\varrho H_N+(1-\varrho)
\frac{N}{N+1}\hat{\mathcal{H}}_N\biggl]
-J'(H_N)\biggl\}d(xf_j(x)) \\
  C_{33N,j}^{**} &=& n_j^{-1/2}\mathcal{A}_j\int_{S_{N_\epsilon}^c}\biggl(\frac{N}{N+1}\mathcal{H}_N-H_N\biggl)
\\&&\times
\biggl\{J'\biggl[\varrho H_N+(1-\varrho)
\frac{N}{N+1}\hat{\mathcal{H}}_N\biggl]
-J'(H_N)\biggl\}d(xf_j(x))
\end{eqnarray*}
Following the arguments of $C_{31N,j}^*$, and
using (A.2)$-$(A.4), we obtain
\begin{eqnarray*}
|C_{33N,j}^*|
&\le&\frac{n_j^{-1/2}}{c_j}|\mathcal{A}_j|
\int_{S_{N_\epsilon}}
\biggl|\frac{N}{N+1}\mathcal{H}_N-H_N
\biggl|\nonumber\\
&&\times\biggl|J'\biggl[\varrho
H_N+(1-\varrho)\frac{N}{N+1
}\hat{\mathcal{H}}_N\biggl]-J'(H_N)\biggl|dF_j(x)\nonumber\\
&\le& \frac{Kn_j^{-1/2}}{c_j}
(1+\zeta^{\delta-\frac32})\mathcal{O}(N^{-1/2})
\int_{\frac{K }{N}}^1 H_N^{\delta-\delta'-1}
dH_N(x).
\end{eqnarray*}
Therefore, $C_{33N,j}^*=o_p(N^{-1/2})$.
Similarly, in view of (22), we can show
$C_{33N,j}^{**}=o_p(N^{-1/2})$. Hence, by (21),
we have $C_{33N,j}=o_p(N^{-1/2})$. To complete
the evaluation of $C_{3N,j}$, we consider
$C_{34N,j}= (C_{34N,j}^*+C_{34N,j}^{**})$, where
\begin{eqnarray*}
  C_{34N,j}^* &=& \frac{n_j^{-1/2}\mathcal{A}_j}{N+1}
\sum_{i=1}^kn_i^{1/2}\int_{S_{N_\epsilon}} x
f_i(x)\\
&&\times\biggl\{J'\biggl[\varrho
H_N+(1-\varrho)\frac{N}{N+1}\hat{\mathcal{H}}_N
\biggl]-J'(H_N)\biggl\}d(xf_j(x)), \\
  C_{34N,j}^{**} &=& \frac{n_j^{-1/2}\mathcal{A}_j}{N+1}
\sum_{i=1}^kn_i^{1/2}\int_{S_{N_\epsilon}^c} x
f_i(x)
\\
&&\times
\biggl\{J'\biggl[\varrho
H_N+(1-\varrho)\frac{N}{N+1}\hat{\mathcal{H}}_N
\biggl]-J'(H_N)\biggl\}d(xf_j(x)).
\end{eqnarray*}
We first turn to evaluate $C_{34N,j}^*$. From
(A.2)$-$(A.4), (21) and (32), it follows that
\begin{eqnarray}
|C_{34N,j}^*| &\le&
K\frac{n_j^{-1/2}|\mathcal{A}_j|}{c_jN}
\sum_{i=1}^kn_i^{1/2}|
\mathcal{A}_i|\int_{S_{N_\epsilon}}
H_N(1-H_N)\nonumber\\
&&\times\biggl|J'\biggl[\varrho
H_N+(1-\varrho)\frac
{N}{N+1}\hat{\mathcal{H}}_N\biggl]-J'(H_N)
\biggl|dF_j(x)\nonumber\\
&\le& \mathcal{O}_p\biggl(n_j^{-1/2}N^{-1}
\sum_{i=1}^kn_i^{1/2}\biggl) \int_{\frac{K}{N}}^1
H_N^{\delta-\frac12}dH_N(x).
\end{eqnarray}
Thus, $C_{34N,j}^*=o_p(N^{-1/2})$. By analogy
with (23) and (33), we can show
$C_{34N,j}^{**}=o_p(N^{-1/2})$. Consequently, we
have
\[
C_{3N,j}=o_p(N^{-1/2}).
\]
This completes the proof of the theorem.
\\
\\
\textbf{Acknowledgement}\\\\
The author is grateful to the associate editor
and two referees for their insightful remarks and
suggestions that greatly improved the original
version of this paper.

\end{document}